\documentclass[12pt]{amsart}
\def\a{\alpha}
\def\b{\beta}

\def\ci{\circ}

\def\D{\Delta}

\def\g{\gamma}

\def\lr{\longrightarrow}
\def\o{\otimes}

\def\s{\sigma}

\def\v{\varepsilon}

\def\1{^{-1}}

\usepackage{amssymb}

\usepackage{enumerate}
\usepackage{mathrsfs}
\usepackage{amssymb,amsfonts}

\usepackage{graphicx}

\makeatletter
\@namedef{subjclassname@2010}{%
  \textup{2010} Mathematics Subject Classification}
\makeatother

\theoremstyle{definition}

\numberwithin{equation}{section}

\begin{document}


\baselineskip=17pt



\title[Cobraided Smash Product Hom-Hopf Algebras]{Cobraided Smash Product Hom-Hopf Algebras}

\author[T. Ma]{Tianshui Ma}
\address{School of Mathematics and Information Science\\ Henan Normal University\\
453007 Xinxiang, China}
\email{matianshui@yahoo.com}

\author[H. Li]{Haiying Li}
\address{School of Mathematics and Information Science\\ Henan Normal University\\
453007 Xinxiang, China}
\email{haiyingli2012@yahoo.com}

\author[T. Yang]{Tao Yang}
\address{College of Science,\\ Nanjing Agricultural University\\
210095 Nanjing, China}
\email{tao.yang.seu@gmail.com}

\date{}

\begin{abstract}
 Let $(A,\a)$ and $(B,\b)$ be two Hom-Hopf algebras. In this paper, we construct a class of new Hom-Hopf
 algebras: $R$-smash product $(A\natural_R B,\a\o \b)$. Moreover, necessary and sufficient conditions for $(A\natural_R B,\a\o \b)$ to
 be a cobraided Hom-Hopf algebra are given.
\end{abstract}

\thanks{Corresponding author: Tianshui Ma}



\maketitle

\section{Introduction}
 Hom-structures (Lie algebras, algebras, coalgebras, Hopf algebras) have been intensively investigated in the literature recently, see \cite{CG,CWZ,MS1,MS2,Yau1,Yau3,Yau4}. Hom-algebras are generalizations of algebras obtained by a twisting map, which have been introduced for the first time in \cite{MS1}. The associativity is replaced by Hom-associativity, Hom-coassociativity for a Hom-coalgebra can be considered in a similar way. Also definitions and properties of Hom-bialgebras and Hom-Hopf algebras have been proposed, see \cite{CG,CWZ,MS2,Yau3,Yau4}.

 In \cite{CG}, Caenepeel and Goyvaerts studied the Hom-structures from the point of view of monoidal categories and introduced that Hom-algebras coincide with algebras in a symmetric monoidal category. In \cite{Yau3}, Yau presented the notion of cobraided Hom-bialgebras
 and got the following conclusion: each cobraided Hom-bialgebra comes with solutions of the operator quantum Hom-Yang-Baxter equations, which are twisted analogues of the operator form of the quantum Yang-Baxter equation. Solutions of the Hom-Yang-Baxter equation can be obtained from comodules of suitable cobraided Hom-bialgebras. And in \cite{Yau1}, Yau introduced and characterized the concept of module Hom-algebras as a twisted version of usual module algebras.

 Let $H$ be a Hopf algebra and $A$ an $H$-module algebra, then, as we all know, we can construct a new Hopf algebra: smash product $A\# H$ (see \cite{Mo} or \cite{Mon}). The extended forms of smash product can be found in the following literature \cite{CIMZ,MW}.

 Let $(H,\b)$ be a Hom-Hopf algebra and $(A,\a)$ an $(H,\b)$-module Hom-algebra (introduced by Yau in \cite{Yau1}), then it is natural to ask: How to construct smash product Hom-Hopf algebra and when the smash product Hom-Hopf algebra is cobraided ?

 The purpose of this article is to give a positive answer to the above questions.

 This article is organized as follows. In Section 2, we recall some definitions and results which will be used later. In Section 3, instead of constructing smash product Hom-Hopf algebra $(A\natural H,\a\o \b)$ (see Theorem 3.3.), we give a more general case, so-called $R$-smash product Hom-Hopf algebra $(A\natural_R B,\a\o \b)$ (see Theorem 3.1).  We remark that the smash product Hom-Hopf algebra $(A\natural H,\a\o \b)$ in Theorem 3.3 is different from the one defined by Chen-Wang-Zhang in \cite{CWZ}, since here the construction of $(A\natural H,\a\o \b)$ is based on the concept of the module Hom-algebra introduced by Yau in \cite{Yau1} which differs from the Chen-Wang-Zhang's in \cite{CWZ}.  Necessary and sufficient conditions for $(A\natural_R B,\a\o \b)$ to be a cobraided Hom-Hopf algebra are derived in Section 4 (see Theorem 4.8, Theorem 4.9). And in last section, we give a concrete example.

\section{Preliminaries}
 Throughout this paper, we follow the definitions and terminologies in \cite{CG,Yau1,Yau3}, with all algebraic systems supposed to be over the field $K$. Given a $K$-space $M$, we write $id_M$ for the identity map on $M$.
\smallskip

 We now recall some useful definitions.

\smallskip

 {\bf Definition 2.1} A Hom-algebra is a quadruple $(A,\mu,1_A,\a)$ (abbr. $(A,\a)$), where $A$ is a $K$-linear space, $\mu: A\o A \lr A$ is a $K$-linear map, $1_A \in A$ and $\a$ is an automorphism of $A$, such that
 \begin{eqnarray*}
 &(HA1)& \a(aa')=\a(a)\a(a');~~\a(1_A)=1_A,\\
 &(HA2)& \a(a)(a'a'')=(aa')\a(a'');~~a1_A=1_Aa=\a(a)
 \end{eqnarray*}
 are satisfied for $a, a', a''\in A$. Here we use the notation $\mu(a\o a')=aa'$.

\smallskip

 {\bf Definition 2.2} A Hom-coalgebra is a quadruple $(C,\D,\v_C,\b)$ (abbr. $(C,\b)$), where $C$ is a $K$-linear space, $\D: C \lr C\o C$, $\v_C: C\lr K$ are $K$-linear maps, and $\b$ is an automorphism of $C$, such that
 \begin{eqnarray*}
 &(HC1)& \b(c)_1\o \b(c)_2=\b(c_1)\o \b(c_2);~~\v_C\ci \b=\v_C\\
 &(HC2)& \b(c_{1})\o c_{21}\o c_{22}=c_{11}\o c_{12}\o \b(c_{2});~~\v_C(c_1)c_2=c_1\v_C(c_2)=\b(c)
 \end{eqnarray*}
 are satisfied for $c \in A$. Here we use the notation $\D(c)=c_1\o c_2$ (summation implicitly understood).

 \smallskip

 {\bf Remarks}  (1) The first equation in $(HC2)$ is equivalent to
 $$
 c_{1}\o c_{21}\o c_{22}= \b^{-1}(c_{11})\o c_{12}\o \b(c_{2}) \eqno(1)
 $$
 and
 $$
 c_{11}\o c_{12}\o c_{2}= \b(c_{1})\o c_{21}\o \b^{-1}(c_{22}), \eqno(2)
 $$
 respectively.\newline
 \indent{\phantom{\bf Remarks}} (2) By $(1),(2)$ and $(HC2)$, we have
 $$
 c_{11}\o c_{12}\o c_{21}\o c_{22}=\b(c_{1})\o \b^{-1}(c_{211})\o \b^{-1}(c_{212})\o c_{22}. \eqno(3)
 $$

\smallskip

 {\bf Definition 2.3} A Hom-bialgebra is a sextuple $(H,\mu,1_H,\D,\v,\g)$ (abbr. $(H,\g)$), where
 $(H,\mu,1_H,\g)$ is a Hom-algebra and $(H,\D,\v,\g)$ is a Hom-coalgebra, such that
 $\D$ and $\v$ are morphisms of Hom-algebras, i.e.
 $$
 \D(hh')=\D(h)\D(h');~~\D(1_H)=1_H\o 1_H,
 $$
 $$
 \v(hh')=\v(h)\v(h');~~\v(1_H)=1.
 $$
 Furthermore, if there exists a linear map $S: H\lr H$ such that
 $$
 S(h_1)h_2=h_1S(h_2)=\v(h)1_H~\hbox{and}~S(\g(h))=\g(S(h)),
 $$
 then we call $(H,\mu,1_H,\D,\v,\g,S)$(abbr. $(H,\g,S)$) a Hom-Hopf algebra.

 Let $(H,\g)$ and $(H',\g')$ be two Hom-bialgebras. The linear map $f: H\lr H'$ is called a Hom-bialgebra map if $f\ci \g=\g'\ci f$
 and at the same time $f$ is a bialgebra map in the usual sense.

 \smallskip

 {\bf Definition 2.4} Let $(A,\b)$ be a Hom-algebra. A left $(A,\b)$-Hom-module is
 a triple $(M,\rhd,\a)$, where $M$ is a linear space, $\rhd: A\o M \lr M$ is a linear map,
 and $\a$ is an automorphism of $M$, such that
 \begin{eqnarray*}
 &(HM1)& \a(a\rhd m)=\b(a)\rhd \a(m),\\
 &(HM2)& \b(a)\rhd (a'\rhd m)=(aa')\rhd \a(m);~~ 1_A\rhd m=\a(m)
 \end{eqnarray*}
 are satisfied for $a, a' \in A$ and $m\in M$.

 \smallskip

 {\bf Remarks} (1) It is obvious that $(A,\mu,\b)$ is a left
 $(A,\b)$-Hom-module.\newline
 \indent{\phantom{\bf Remarks}} (2) When $\b=id_A$ and $\a=id_M$, a left $(A,\b)$-Hom-module
 is the usual left $A$-module.

 \smallskip

 {\bf Definition 2.5} Let $(H,\b)$ be a Hom-bialgebra and $(A,\a)$ a Hom-algebra. If
 $(A,\rhd,\a)$ is a left $(H,\b)$-Hom-module and for all $h\in H$ and $a, a'\in A$,
 \begin{eqnarray*}
 &(HMA1)& \b^{2}(h)\rhd (aa')=(h_1\rhd a)(h_2\rhd a'),\\
 &(HMA2)& h\rhd 1_A=\v_H(h)1_A,
 \end{eqnarray*}
 then $(A,\rhd,\a)$ is called an $(H,\b)$-module Hom-algebra.

 \smallskip

 {\bf Remarks} (1) When $\a=id_A$ and $\b=id_H$, an $(H,\b)$-module Hom-algebra
 is the usual $H$-module algebra.\newline
 \indent{\phantom{\bf Remarks}} (2) Similar to the case of Hopf
 algebras, in \cite{Yau4}, Yau concluded that the Eq.$(HMA1)$ is satisfied if and only if $\mu_A$ is a
 morphism of $H$-modules for suitable $H$-module structures on $A\o A$ and $A$, respectively.\newline
 \indent{\phantom{\bf Remarks}} (3) If $\b^2=id$ in $(HMA1)$, then we can get $(6.1)$ in \cite{CWZ}. So two definitions of module Hom-algebra are different which leads to the differentness of smash product Hom-Hopf algebra in Theorem 3.3 and in Definition 6.2 in \cite{CWZ}.

 \smallskip

 {\bf Definition 2.6} A cobraided Hom-Hopf algebra is a octuple $(H,\mu,1_H, \D,$ $\v,S,\a,\s)$
 (abbr.$(H,\a,\s)$) in which $(H,\mu,1_H,\D,\v,S,\a)$ is a Hom-Hopf algebra and $\s$ is a bilinear form on $H$
 (i.e., $\s \in Hom(H\o H, K)$), satisfying the following axioms (for all $h, g, l\in H$):
 \begin{eqnarray*}
 &(CHA1)& \s(h,1_H)=\s(1_H,h)=\v(h),\\
 &(CHA2)& \s(hg,\a(l))=\s(\a(h),l_1)\s(\a(g),l_2),\\
 &(CHA3)& \s(\a(h),gl)=\s(h_1,\a(l))\s(h_2,\a(g)),\\
 &(CHA4)& \s(h_1,g_1)h_2g_2=g_1h_1\s(h_2,g_2),\\
 &(CHA5)& \s(\a(h),\a(g))=\s(h,g).
 \end{eqnarray*}
 In this case, $\s$ is called the Hom-cobraiding form.

 \smallskip

 {\bf Remarks} (1) When $\a=id_H$, a cobraided Hom-Hopf algebra
 is exactly the usual cobraided (or coquasitriangular) Hopf algebra.\newline
 \indent{\phantom{\bf Remarks}} (2) It is slightly different from
 the definition in \cite{Yau3} or \cite{Yau4}. Here we replace the
 Hom-bialgebra by Hom-Hopf algebra and also add another two conditions
 $(CHA1)$ and $(CHA5)$. Similar to the Hopf algebra setting, the Hom-cobraiding form
 $\s$ in Definition 2.6 is invertible.\newline
 \indent{\phantom{\bf Remarks}} (3) Based on Yau's results in
 \cite{Yau3}, each cobraided Hom-Hopf algebra comes with solutions of the operator
 quantum Hom-Yang-Baxter equations, which are twisted analogues of the operator form
 of the quantum Yang-Baxter equation.
\smallskip

 Next, we generalize the concept of skew pairing to Hom-setting.

 {\bf Definition 2.7}  Let $(A,\a,S_A)$ and $(B,\b,S_B)$ be two Hom-Hopf algebras,
 $\vartheta\in Hom(A\o B, K)$ a bilinear form. A Hom-skew pairing is a triple $(A,B,\vartheta)$
 such that
 \begin{eqnarray*}
 &(SP1)& \vartheta(a,1_B)=\v_A(a);~~\vartheta(1_A,b)=\v_B(b),\\
 &(SP2)& \vartheta(aa',\b(b))=\vartheta(\a(a),b_1)\vartheta(\a(a'),b_2),\\
 &(SP3)& \vartheta(\a(a),bb')=\vartheta(a_1,\b(b'))\vartheta(a_2,\b(b)),\\
 &(SP4)& \vartheta(\a(a),\b(b))=\vartheta(a,b),
 \end{eqnarray*}
 where $a,a' \in A$ and $b,b'\in B$.

 \smallskip

 {\bf Remarks} (1) When $\a=id_A$ and $\b=id_B$, we can get the usual skew pairing.\newline
 \indent{\phantom{\bf Remarks}} (2) If $(H,\a,\s)$ is a cobraided Hom-Hopf algebra,
 then $(H,H,\s)$ is a Hom-skew pairing.\newline
 \indent{\phantom{\bf Remarks}} (3) $\vartheta$ is (convolution) invertible with
 $\vartheta\1(a,b)=\vartheta(S_A(a),b)$.

\section{Smash product Hom-Hopf algebra}

 In this section, we introduce a class of new Hom-Hopf algebras: $R$-smash product $A\natural_R B$, generalizing $R$-smash product studied in \cite{CIMZ}. As a special case, Hom-smash product is derived based on the structure of module Hom-algebra introduced by Yau in \cite{Yau1} or \cite{Yau4}.

 Let $A$ and $B$ be two linear spaces, $R: B\o A \lr A\o B$ a linear map. In the following, we write $R(b\o a)=\sum a_R\o b_R$ for all $a\in A$ and $b\in B$, and the notations $\sum a_r\o b_r$, $\sum a_{R'}\o b_{R'}$ are the copies of $\sum a_R\o b_R$.  As usual, we omit the summation sign ``$\sum$".

 {\bf Theorem 3.1} Let $(A,\mu_A,1_A,\a)$ and $(B,\mu_B,1_B,\b)$ be
 two Hom-algebras, $R: B\o A \lr A\o B$ a linear map such that for all $a\in A, b\in B$,
 $$
 \a(a)_R\o \b(b)_R=\a(a_R)\o \b(b_R). \eqno(4)
 $$
 Then $(A\natural_R B, \a\o \b)$ ($A\natural_R B=A\o B$ as a linear space) with the multiplication
 $$
 (a\o b)(a'\o b')=a\a^{-1}(a')_{R}\o \b^{-1}(b_{R})b',
 $$
 where $a, a'\in A, b, b'\in B$, and unit $1_A\o 1_B$ becomes a Hom-algebra if and only if the
 following conditions hold:
 \begin{eqnarray*}
 &(C1)& a_{R}\o 1_{BR}=\a(a)\o 1_B;~~1_{AR}\o b_R=1_A\o \b(b),\\
 &(C2)& \a(a)_R\o (bb')_R=a_{Rr}\o \b^{-1}(\b(b)_r)b'_R,\\
 &(C3)& \a((aa')_R)\o \b(b)_R=\a(a_R)\a(a')_r\o b_{Rr},
 \end{eqnarray*}
 where $a, a'\in A, b, b'\in B$.

 We call this Hom-algebra $R$-smash product Hom-algebra and denote it by $(A\natural_R B,\a\o \b)$.

 {\bf Proof} ($\Longleftarrow$) For all $a, a', a''\in A$, $b, b', b''\in
 B$, firstly, we prove that $(HA1)$ holds. In fact, one
 can get
 \begin{eqnarray*}
 (\a\o \b)((a\o b)(a'\o b'))
 &\stackrel{}{=}&\underline{\a(a\a\1(a')_{R})}\o \underline{\b(\b\1(b_{R})b')}\\
 &\stackrel{(HA1)}{=}&\a(a)\underline{\a(\a\1(a')_{R})}\o \underline{b_{R}}\b(b')\\
 &\stackrel{(4)}{=}&\a(a)a'_{R}\o \b\1(\b(b)_{R})\b(b')\\
 &\stackrel{}{=}&((\a\o \b)(a\o b))((\a\o \b)(a'\o b'))
 \end{eqnarray*}
 and
 $$
 (\a\o \b)(1_A\o 1_B)=\a(1_A)\o \b(1_B)\stackrel{(HA1)}{=}1_A\o 1_B.
 $$
 Secondly, we compute the condition $(HA2)$ as follows.
 \begin{eqnarray*}
 &&(\a(a)\o \b(b))((a'\o b')(a''\o b''))\\
 &\stackrel{}{=}&\a(a)\underline{\a\1(a'\a\1(a'')_{R})_{r}}
 \o \b\1(\underline{\b(b)_{r}})(\b\1(b'_{R})b'')\\
 &\stackrel{(C3)}{=}&\a(a)\a\1(\a(\a\1(a')_{r})\a\1(a'')_{RR'})
 \o \b\1(b_{rR'})(\b\1(b'_{R})b'')\\
 &\stackrel{}{=}&\underline{\a(a)(\a\1(a')_{r}\a\1(\a\1(a'')_{RR'}))}
 \o \b\1(b_{rR'})(\b\1(b'_{R})b'')\\
 &\stackrel{(HA2)}{=}&(a\a\1(a')_{r})\a\1(a'')_{RR'}
 \o \underline{\b\1(b_{rR'})(\b\1(b'_{R})b'')}\\
 &\stackrel{(HA2)}{=}&(a\a\1(a')_{r})\underline{\a\1(a'')_{RR'}}
 \o \b\1(\underline{\b\1(b_{rR'})b'_{R}})\b(b'')\\
 &\stackrel{(C2)}{=}&(a\a\1(a')_{r})a''_{R}\o \b\1((\b\1(b_{r})b')_{R})\b(b'')\\
 &\stackrel{}{=}&((a\o b)(a'\o b'))(\a(a'')\o \b(b''))
 \end{eqnarray*}
 and
 \begin{eqnarray*}
 (a\o b)(1_A\o 1_B)
 &\stackrel{}{=}&a\a^{-1}(1_A)_{R}\o \b^{-1}(b_{R})1_B\\
 &\stackrel{(HA1)}{=}&a1_{AR}\o \b^{-1}(b_{R})1_B\\
 &\stackrel{(C1)}{=}&a1_{A}\o b1_B\\
 &\stackrel{(HA2)}{=}&\a(a)\o \b(b).
 \end{eqnarray*}
 Similarly, $(1_A\o 1_B)(a\o b)=\a(a)\o \b(b)$ holds.

 ($\Longrightarrow$) By $(HA2)$, we have
 $$
 1_A\a\1(a)_R\o \b\1(1_{BR})b=\a(a)\o \b(b), \eqno(5)
 $$
 $$
 a\a\1(1_A)\o \b\1(b_R)1_B=\a(a)\o \b(b) \eqno(6)
 $$
 and
 \begin{eqnarray*}
 &&\a(a)\a\1(a'\a\1(a'')_{R})_{r}\o \b\1(\b(b)_{r})(\b\1(b'_{R})b'')\\
 &&~~~~=(a\a\1(a')_{r})a''_{R}\o \b\1((\b\1(b_{r})b')_{R})\b(b'').\ \ \ \ \ \ \ \ \ \ \ \ \ \ \ \ \ \ \ \ \ \ \ \ \ \ \ \ \ \ \ \ \ (7)
 \end{eqnarray*}

 Let $b=1_B$ and $a=1$ in Eqs.$(5)$ and $(6)$, respectively,
 we can get $(C1)$.

 Let $a=a'=1_A$ and $b''=1_B$ in Eq.$(7)$ and by $(C1)$, then
 $(C2)$ holds.

 Likewise, $(C3)$ can be obtained by letting $a=1_A$ and $b'=b''=1_B$ in Eq.$(7)$. \hfill $\square$

\smallskip

 When $\a=id_A$ and $\b=id_B$, we have

 {\bf Example 3.2}(\cite{CIMZ})  Let $(A,\mu_A,1_A)$ and $(B,\mu_B,1_B)$ be two algebras, $R: B\o A \lr A\o B$ a linear map. Then $A\#_R B$ ($A\#_R B=A\o B$ as a linear
 space) with the multiplication
 $$
 (a\o b)(a'\o b')=aa'_{R}\o b_{R}b',
 $$
 where $a, a'\in A, b, b'\in B$, and unit $1_A\o 1_B$ becomes an algebra if and only if the
 following conditions hold:
 \begin{eqnarray*}
 &(1)& a_{R}\o 1_{BR}=a\o 1_B;~~1_{AR}\o b_R=1_A\o b,\\
 &(2)& a_R\o (bb')_R=a_{Rr}\o b_rb'_R,\\
 &(3)& (aa')_R\o b_R=a_Ra'_r\o b_{Rr},
 \end{eqnarray*}
 where $a, a'\in A, b, b'\in B$.

\smallskip

 {\bf Theorem 3.3} Let $(H,\b)$ be a Hom-bialgebra and $(A,\rhd,\a)$ an $(H,\b)$-module
 Hom-algebra. Then $(A\natural H, \a\o \b)$ ($A\natural H=A\o H$ as a linear space) with the multiplication
 $$
 (a\o h)(a'\o h')=a(h_1\rhd \a^{-1}(a'))\o \b^{-1}(h_{2})h',
 $$
 where $a, a'\in A, h, h'\in H$, and unit $1_A\o 1_H$ is a
 Hom-algebra, we call it smash product Hom-algebra denoted by $(A\natural H,\a\o \b)$.

 {\bf Proof} Define $R:H\o A\lr A\o H$ by
 $$
 R(h\o a)=h_{1}\rhd a\o h_{2}, \forall a\in A, h\in H.
 $$

 Firstly, for all $a\in A$ and $h\in H$,
 \begin{eqnarray*}
 \a(a)_R\o \b(h)_R
 &\stackrel{}{=}&\underline{\b(h)_1}\rhd \a(a)\o \underline{\b(h)_2}\\
 &\stackrel{(HC1)}{=}&\underline{\b(h_1)\rhd \a(a)}\o \b(h_2)\\
 &\stackrel{(HM1)}{=}&\a(h_1\rhd a)\o \b(h_2)=\a(a_R)\o \b(h_R),
 \end{eqnarray*}
 so Eq.$(4)$ holds.

 Secondly, we have
 $$
 a_{R}\o 1_{HR}=1\rhd a\o 1_H\stackrel{(HM2)}{=}\a(a)\o 1_H
 $$
 and
 $$
 1_{AR}\o h_R=h_1\rhd 1_A\o h_2\stackrel{(HMA2)}{=}1_A\o
 \v(h_1)h_2\stackrel{(HC2)}{=}1_A\o \b(h).
 $$

 Thirdly, we verify that the conditions $(C2)$ and $(C3)$ are satisfied.
 For all $a, a'\in A, h, h'\in B$,
 \begin{eqnarray*}
 \a(a)_R\o (hh')_R
 &\stackrel{}{=}&(hh')_1\rhd \a(a)\o (hh')_2\\
 &\stackrel{}{=}&\underline{(h_1h'_1)\rhd \a(a)}\o h_2h'_2\\
 &\stackrel{(HM2)}{=}&\underline{\b(h_1)}\rhd (h'_1\rhd a)\o \underline{h_2}h'_2\\
 &\stackrel{(HC1)}{=}&\b(h)_1\rhd (h'_1\rhd a)\o \b^{-1}(\b(h)_2)h'_2\\
 &\stackrel{}{=}&a_{Rr}\o \b^{-1}(\b(h)_r)h'_R
 \end{eqnarray*}
 and
 \begin{eqnarray*}
 \a((aa')_R)\o \b(h)_R
 &\stackrel{}{=}&\a(\underline{\b(h)_1}\rhd (aa'))\o \underline{\b(h)_2}\\
 &\stackrel{(HC1)}{=}&\underline{\a(\b(h_1)\rhd (aa'))}\o \b(h_2)\\
 &\stackrel{(HM1)}{=}&\b^2(h_1)\rhd \underline{\a(aa')}\o \b(h_2)\\
 &\stackrel{(HA1)}{=}&\underline{\b^2(h_1)\rhd (\a(a)\a(a'))}\o \b(h_2)\\
 &\stackrel{(HMA1)}{=}&(\underline{h_{11}}\rhd \a(a))(\underline{h_{12}}\rhd \a(a'))\o \underline{\b(h_2)}\\
 &\stackrel{(HC2)}{=}&\underline{(\b(h_{1})\rhd \a(a))}(h_{21}\rhd \a(a'))\o h_{22}\\
 &\stackrel{(HM1)}{=}&\a(h_{1}\rhd a)(h_{21}\rhd \a(a'))\o h_{22}\\
 &\stackrel{}{=}&\a(a_R)\a(a')_r\o h_{Rr}.
 \end{eqnarray*}
 Thus we complete the proof. \hfill $\square$

 \smallskip

 {\bf Remarks} (1) The smash product Hom-Hopf algebra $(A\natural H,\a\o \b)$ is different from the one defined by Chen-Wang-Zhang in \cite{CWZ}, since here the construction of $(A\natural B,\a\o \b)$ is based on the concept of the module Hom-algebra introduced by Yau in \cite{Yau1}, while two of conditions $(6.1),(6.2)$ in the module Hom-algebra in \cite{CWZ} are same to the case of Hopf algebra.  \newline
 \indent{\phantom{\bf Remarks}} (2) When $\a=id_A$ and $\b=id_H$, we can get the usual smash product algebra $A\# H$ (see \cite{Mo,Mon}).

 \smallskip

 {\bf Lemma 3.4} Let $(C,\a)$ and $(D,\b)$ be two Hom-coalgebras.
 Then $(C\o D, \a\o \b)$ is a Hom-coalgebra with the following
 comultiplication and counit
 $$
 \D(c\o d)=c_{1}\o d_{1}\o c_{2}\o d_{2},
 $$
 $$
 \v(c\o d)=\v_C(c)\v_D(d),
 $$
 in which $c\in C$ and $d\in D$. We call it tensor product
 Hom-coalgebra.

 {\bf Proof} Straightforward. \hfill $\square$

\smallskip

 {\bf Theorem 3.5} Let $(A,\a,S_A)$ and $(B,\b,S_B)$ be two Hom-Hopf algebras, $R: B\o A \lr A\o B$ a linear map.
 Then the $R$-smash product Hom-algebra $(A\natural_R B,\a\o \b)$ equipped with the tensor product Hom-coalgebra structure becomes a Hom-bialgebra if and only if $R$ is a coalgebra map, i.e.
 $$
 a_{R1}\o b_{R1}\o a_{R2}\o b_{R2}=a_{1R}\o b_{1R}\o a_{2r}\o b_{2r},
 $$
 $$
 \v_A(a_R)\v_B(b_R)=\v_A(a)\v_B(b),
 $$
 where $a\in A$, $b\in B$.

 Furthermore, $R$-smash product Hom-bialgebra $(A\natural_R B,\a\o \b)$ is a Hom-Hopf algebra with antipode $\bar{S}$ defined by
 $$
 \bar{S}(a\o b)=\a\1(S_A(a))_R\o \b\1(S_B(b)_R).
 $$

 {\bf Proof} We only prove that $\bar{S}$ is an antipode of $(A\natural_R B,\a\o \b)$.
 The rest is straightforward by direct computation. For all $a\in A$ and $b\in B$,
 \begin{eqnarray*}
 (\bar{S}*id_{A\natural_R B})(a\o b)
 &\stackrel{}{=}&(\a\1(S_A(a_1))_R\o \b\1(S_B(b_1)_R))(a_2\o b_2)\\
 &\stackrel{}{=}&\a\1(S_A(a_1))_R\underline{\a\1(a_2)_r}\o \b\1(\underline{\b\1(S_B(b_1)_R)_r})b_2\\
 &\stackrel{(4)}{=}&\underline{\a\1(S_A(a_1))_R\a\1(a_{2r})}\o \b^{-2}(S_B(b_1)_{Rr})b_2\\
 &\stackrel{(HA1)}{=}&\a\1(\underline{\a(\a\1(S_A(a_1))_R)a_{2r}})\o \b^{-2}(\underline{S_B(b_1)_{Rr}})b_2\\
 &\stackrel{(C3)}{=}&\a\1(S_A(a_1)a_{2})_R\o \b^{-2}(\b(S_B(b_1))_{R})b_2\\
 &\stackrel{}{=}&\underline{1_{AR}}\v_A(a)\o \b^{-2}(\underline{\b(S_B(b_1))_{R}})b_2\\
 &\stackrel{(C1)}{=}&1_{A}\v_A(a)\o S_B(b_1)b_2\\
 &\stackrel{}{=}&1_{A}\o 1_{B}\v_A(a)\v_B(b)\\
 &\stackrel{}{=}&1_{A}\o 1_{B}\bar{\v}(a\o b)
 \end{eqnarray*}
 and
 \begin{eqnarray*}
 (id_{A\natural_R B}*\bar{S})(a\o b)
 &\stackrel{}{=}&(a_1\o b_1)(\a\1(S_A(a_2))_R\o \b\1(S_B(b_2)_R))\\
 &\stackrel{}{=}&a_1\underline{\a\1(\a\1(S_A(a_2))_{R})}_{r}\o \b\1(b_{1r})\underline{\b\1(S_B(b_2)_R)}\\
 &\stackrel{(4)}{=}&a_1\a^{-2}(S_A(a_2))_{Rr}\o \b\1(b_{1r})\b\1(S_B(b_2))_R\\
 &\stackrel{}{=}&a_1\underline{\a\1(\a\1(S_A(a_2)))_{Rr}}\\
 && \o \underline{\b\1(\b(\b\1(b_1))_{r})\b\1(S_B(b_2))_R}\\
 &\stackrel{(C2)}{=}&a_1\a\1(S_A(a_2))_{R}\o \b\1(b_1S_B(b_2))\\
 &\stackrel{}{=}&a_1\underline{\a\1(S_A(a_2))_{R}}\o \underline{1_{BR}}\v_B(b)\\
 &\stackrel{(C1)}{=}&a_1S_A(a_2)\o 1_{B}\v_B(b)\\
 &\stackrel{}{=}&1_{A}\o 1_{B}\v_A(a)\v_B(b)\\
 &\stackrel{}{=}&1_{A}\o 1_{B}\bar{\v}(a\o b),
 \end{eqnarray*}
 while

  \begin{eqnarray*}
 \bar{S}(\a(a)\o \b(b))
 &\stackrel{}{=}&\a\1(S_A(\a(a)))_R\o \b\1(S_B(\b(b))_R)\\
 &\stackrel{}{=}&\a\1(\a(S_A(a)))_R\o \b\1(\b(S_B(b))_R)\\
 &\stackrel{}{=}&S_A(a)_R\o \b\1(\b(S_B(b))_R)\\
 &\stackrel{(4)}{=}&\a(\a\1(S_A(a))_R)\o S_B(b)_R\\
 &\stackrel{}{=}&(\a\o \b)(\bar{S}(a\o b)),
 \end{eqnarray*}
 finishing the proof.   \hfill $\square$

 \smallskip
 When $\a=id_A$ and $\b=id_B$, we have

 {\bf Example 3.6}(\cite{CIMZ}) Let $A$ and $B$ be two Hopf algebras.
 Then the twisted tensor product algebra $A\#_R B$ equipped with
 the usual tensor product coalgebra structure is a bialgebra if and only if
 $R$ is a coalgebra map.

 Furthermore, twisted tensor product bialgebra $A\#_R B$ is a
 Hopf algebra with antipode $S_{A\#_R B}$ defined by
 $$
 S_{A\#_R B}(a\o b)=S_A(a)_R\o S_B(b)_R.
 $$

 \smallskip

 {\bf Theorem 3.7} Let $(H,\b)$ be a Hom-Hopf algebra and $(A,\rhd,\a)$ an $(H,\b)$-module
 Hom-algebra. Then the smash product Hom-algebra $(A\natural H, \a\o \b)$ endowed with the tensor product Hom-coalgebra structure
 becomes a Hom-bialgebra if and only if
 $$
 (h\rhd a)_1\o (h\rhd a)_2=(h_1\rhd a_1)\o (h_2\rhd a_2);~~\v_A(h\rhd
 a)=\v_A(a)\v_H(h) \eqno(8)
 $$
 and
 $$
 h_1\o h_2\rhd a=h_2\o h_1\rhd a. \eqno(9)
 $$

 Moreover, smash product Hom-bialgebra $(A\natural H,\a\o \b)$ is a Hom-Hopf algebra with antipode
 $$
 S_{A\natural H}(a\o h)=S_H(h)_1\rhd \a\1(S_A(a))\o \b\1(S_H(h)_2).
 $$

 {\bf Proof} Let $R(h\o a)=h_{1}\rhd a\o h_{2}, \forall a\in A, h\in H$ in Theorem 3.5.
 Then $R$ is a coalgebra map if and only if
 $$
 (h_{1}\rhd a)_1\o h_{21}\o (h_{1}\rhd a)_2\o h_{22}
 =h_{11}\rhd a_1\o h_{12}\o h_{21}\rhd a_2\o h_{22} \eqno(10)
 $$
 and
 $$
 \v_A(h\rhd a)=\v_A(a)\v_{H}(h).
 $$
 And by Eq.$(3)$ and $(HC1)$, it is easy to obtain that the first
 equation in Eq.$(8)$ and Eq.$(9)$ are equivalent to Eq.(10). \hfill $\square$

\smallskip

 {\bf Remarks} (1) Let $(H,\b)$ be a Hom-Hopf algebra. Assume that $(A,\rhd,\a)$ is a Hom-coalgebra
 and an $(H,\b)$-Hom-module satisfying Eq.$(8)$. Then we call $(A,\rhd,\a)$ an $(H,\b)$-module Hom-coalgebra.

 When $\a=id_A$ and $\b=id_H$, then an $(H,\b)$-module Hom-coalgebra is exactly the
 module coalgebra in the usual case (see \cite{Mo}). \newline
 \indent{\phantom{\bf Remarks}} (2) Theorem 3.7 is the Hom-version
 of the usual smash product Hopf algebra (see \cite{Mo}).
\smallskip

\section{Cobraided Hom-Hopf algebra}
\def\theequation{4. \arabic{equation}}
\setcounter{equation} {0} \hskip\parindent
 In this section, necessary and sufficient conditions for smash product Hom-Hopf algebra to be cobraided are given.

 {\bf Proposition 4.1} Let $(A\natural_R B,\a\o \b)$ be a $R$-smash product Hom-Hopf algebra.
 Define
 $$
 i: A\lr A\natural_R B,~ i(a)=a\o 1_B;~~~j: B\lr A\natural_R B, ~j(b)=1_A\o b
 $$
 for all $a\in A$ and $b\in B$. Then $i$ and $j$ are both Hom-bialgebra maps.

 {\bf Proof} Straightforward. \hfill $\square$

\smallskip

 Let $(A\natural_R B,\a\o \b)$ be a $R$-smash product Hom-Hopf algebra, and
 $\s: A\natural_R B\o A\natural_R B\lr K$ a bilinear form. Define
 \begin{eqnarray*}
 &&\tau: A\o A\lr K,~ \tau(a, a')=\s(i\o i)(a\o a'),\\
 &&\upsilon: B\o B\lr K,~ \upsilon(b, b')=\s(j\o j)(b\o b'),\\
 &&\varphi: A\o B\lr K,~ \varphi(a, b)=\s(i\o j)(a\o b),\\
 &&\psi: B\o A\lr K,~ \psi(b, a)=\s(j\o i)(b\o a),
 \end{eqnarray*}
 where $a,a'\in A$ and $b,b'\in B$.

\smallskip

 The following two lemmas are obvious.

 {\bf Lemma 4.2} Let $(A\natural_R B,\a\o \b)$ be a $R$-smash product Hom-Hopf algebra.
 If $\s$ satisfies $(CHA1)$, then for $a\in A$ and $b\in B$,
 \begin{eqnarray*}
 &&(1)\ \ \tau(1_A, a)=\tau(a, 1_A)=\v_{A}(a),\\
 &&(2)\ \ \upsilon(b, 1_B)=\upsilon(1_B, b)=\v_{B}(b), \\
 &&(3)\ \ \varphi(1_A, b)=\v_{B}(b);~\varphi(a, 1_B)=\v_{A}(a),\\
 &&(4)\ \ \psi(1_B, a)=\v_{A}(a);~\psi(b, 1_A)=\v_{B}(b).
 \end{eqnarray*}

 \smallskip

 {\bf Lemma 4.3} Let $(A\natural_R B,\a\o \b)$ be a $R$-smash product Hom-Hopf algebra.
 If $\s$ satisfies $(CHA5)$ for $\a\o \b$, then, for $a,a'\in A$ and $b,b'\in
 B$,
 \begin{eqnarray*}
 &&(1)\ \ \tau(\a(a), \a(a'))=\tau(a, a'),\\
 &&(2)\ \ \upsilon(\b(b), \b(b'))=\upsilon(b, b'),\\
 &&(3)\ \ \varphi(\a(a), \b(b))=\varphi(a, b),\\
 &&(4)\ \ \psi(\b(b), \a(a))=\psi(b, a).
 \end{eqnarray*}

 \smallskip

 {\bf Lemma 4.4} Let $(A\natural_R B,\a\o \b,\s)$ be a cobraided $R$-smash product
 Hom-Hopf algebra. Then, for all $a,a'\in A$ and $b,b'\in B$,
 $$
 \s(\a(a)\o \b(b), \a(a')\o \b(b'))
 =\varphi(a_1, b'_1)\tau(a_2, a'_1)\upsilon(b_1, b'_2)\psi(b_2, a'_2). \eqno(11)
 $$

 {\bf Proof} By $(CHA2)$ and $(CHA3)$, for all $a,a',a'',a'''\in A$,
 $b,b',b'',b'''\in B$, we have
 \begin{eqnarray*}
 &&~~\s(a\a\1(a')_R\o \b\1(b_R)b', a''\a\1(a''')_r\o
 \b\1(b''_r)b''')\\
 &&=\s(a_1\o b_1, a'''_1\o b'''_1)\s(a_2\o b_2, a''_1\o b''_1)\\
 &&\ \ \ \ \times \s(a'_1\o b'_1, a'''_2\o b'''_2)\s(a'_2\o b'_2, a''_2\o b''_2).
 \end{eqnarray*}

 Let $a'=a'''=1_A$ and $b=b''=1_B$ in the above equation, then we
 can get $(11)$.  \hfill $\square$

 \smallskip

 {\bf Lemma 4.5} Let $(A\natural_R B,\a\o \b,\s)$ be a cobraided $R$-smash product
 Hom-Hopf algebra. Then, for all $a,a'\in A$ and $b,b'\in B$,
 \begin{eqnarray*}
 &(D1)& \varphi(\a(\a\1(a)_R), b_1)\upsilon(b'_R, b_2)=\upsilon(\b(b'), b_1)\varphi(\a(a), b_2),\\
 &(D2)& \tau(\a(\a\1(a)_R), a'_1)\psi(b_R, a'_2)=\psi(\b(b), a'_1)\tau(\a(a), a'_2),\\
 &(D3)& \upsilon(b_1, b'_R)\psi(b_2, \a(\a\1(a)_R))=\psi(b_1, \a(a)\upsilon(b_2, \b(b')),\\
 &(D4)& \varphi(a_1, b_R)\tau(a_2, \a(\a\1(a')_R))=\tau(a_1, \a(a')\varphi(a_2, \b(b)),\\
 &(D5)& \psi(b_1, a_1)(\a(\a\1(a_2)_R)\o b_{2R})=(\a(a_1)\o \b(b_{1}))\psi(b_2, a_2),\\
 &(D6)& \varphi(a_1, b_1)(\a(a_2)\o \b(b_{2}))=(\a(\a\1(a_1)_R)\o b_{1R})\varphi(a_2, b_2).
 \end{eqnarray*}

 {\bf Proof} By $(CHA2)$, for all $a,a',a''\in A$,
 $b,b',b''\in B$, we can obtain
 \begin{eqnarray*}
 &&\s(a\a\1(a')_R\o \b\1(b_R)b', \a(a'')\o \b(b''))\\
 &&=\s(\a(a)\o \b(b), a''_1\o b''_1)\s(\a(a')\o \b(b'), a''_2\o
 b''_2). \ \ \ \ \ \ \ \ \ \ \ \ \ \ \ \ \ \ \ \ \ \ (12)
 \end{eqnarray*}
 Let $a=1_A$ and $b'=b''=1_B$ in Eq.$(12)$, then $(C1)$ holds by Eq.$(11)$. Similarly,
 setting $a=a''=1_A$ and $b'=1_B$ in Eq.$(12)$, then we can get $(C2)$ by Eq.$(11)$.

 By $(CHA3)$, for all $a,a',a''\in A$, $b,b',b''\in B$, we have
 \begin{eqnarray*}
 &&\s(\a(a)\o \b(b), a'\a\1(a'')_R\o \b\1(b'_R)b'')\\
 &&=\s(a_1\o b_1, \a(a'')\o \b(b''))\s(a_2\o b_2, \a(a')\o \b(b')).\ \ \ \ \ \ \ \ \ \ \ \ \ \ \ \ \ \ \ \ (13)
 \end{eqnarray*}
 $(C3)$ can be obtained by letting $a=a'=1_A$ and $b''=1_B$ in Eq.$(13)$ and by Eq.$(11)$.
 Likewise, one gets $(C4)$ by putting $a'=1_A$ and $b=b''=1_B$ in Eq.$(13)$ and by Eq.$(11)$.

 By $(CHA4)$, for all $a,a'\in A$, $b,b'\in B$, we have
 \begin{eqnarray*}
 &&\s(a_1\o b_1, a'_1\o b'_1)(a_2\a\1(a'_2)_R\o \b\1(b_{2R})b'_2)\\
 &&=(a'_1\a\1(a_1)_R\o \b\1(b'_{1R})b_1)\s(a_2\o b_2, a'_2\o b'_2).\ \ \ \ \ \ \ \ \ \ \ \ \ \ \ \ \ \ \ \ \ \ \ \ \ \ (14)
 \end{eqnarray*}
 Let $a=1_A$ and $b'=1_B$ in Eq.$(14)$, then we get $(C5)$. And
 $(C6)$ is derived by letting $a'=1_A$ and $b=1_B$ in Eq.$(14)$. \hfill $\square$

 {\bf Lemma 4.6} Given the cobraiding $\s$ on a $R$-smash product Hom-Hopf algebra $(A\natural_R B,\a\o \b)$, consider the induced maps $\tau, \upsilon, \varphi$ and $\psi$. Then

 (1) $(A,\a,\tau)$ and $(B,\b,\upsilon)$ are cobraided Hom-Hopf algebras,

 (2) $(A,B,\varphi)$ and $(B,A,\psi)$ are Hom-skew pairings.

 {\bf Proof} (1) Set $b=b'=b''=1_B$ in Eq.$(12)$ and Eq.$(13)$, we can get
 $(CHA2)$ and $(CHA3)$ for $\tau$, respectively. $(CHA4)$ can be derived by letting $b=b'=1_B$ in Eq.$(14)$,
 then by Lemma 4.2 and Lemma 4.3, $(A,\a,\tau)$ is a cobraided Hom-Hopf
 algebra. Similarly, we can prove that $(B,\b,\upsilon)$ is a cobraided Hom-Hopf algebra.

 (2) Let $a''=1_A$ and $b=b'=1_B$ in Eq.$(12)$, $a'=a''=1_A$ and $b=1_B$ in
 Eq.$(13)$, $(SP2)$ and $(SP3)$ can be obtained for $\varphi$,
 respectively. Then $(A,B,\varphi)$ is a Hom-skew pairing by Lemma 4.2 and Lemma
 4.3. The rest of (2) can be similarly demonstrated.  \hfill $\square$

 {\bf Lemma 4.7} Let $(A\natural_R B,\a\o \b)$ be a $R$-smash product Hom-Hopf
 algebra. If there exist forms $\tau: A\o A\lr K$, $\varphi: A\o B\lr K$, $\psi: B\o A\lr K$,
 and $\upsilon: B\o B\lr K$ such that

 (1) $(A,\a,\tau)$ and $(B,\b,\upsilon)$ are cobraided Hom-Hopf algebras,

 (2) $(A,B,\varphi)$ and $(B,A,\psi)$ are Hom-skew pairings,

 (3) The conditions $(D1)-(D6)$ in Lemma 4.5 hold.\\
 Then $(A\natural_R B,\a\o \b,\s)$ is a cobraided Hom-Hopf algebra with the cobraided structure
 given by
 $$
 \s(\a(a)\o \b(b), \a(a')\o \b(b'))=\varphi(a_1, b'_1)\tau(a_2, a'_1)\upsilon(b_1, b'_2)\psi(b_2,
 a'_2),
 $$
 for $a,a'\in A$ and $b,b'\in B$.

 {\bf Proof} It is obvious that $\s$ satisfies $(CHA1)$ and $(CHA5)$.

 Next, we show that $(CHA2)$ holds for $\s$. For all $a,a',a''\in A$,
 $b,b',b''\in B$,
 \begin{eqnarray*}
 &&\s((a\o b)(a'\o b'), \a(a'')\o \b(b''))\\
 &\stackrel{}{=}&\s(a\a\1(a')_R\o \b\1(b_R)b', \a(a'')\o \b(b''))\\
 &\stackrel{}{=}&\varphi(\underline{\a\1(a\a\1(a')_{R})_{1}}, b''_{1})
 \tau(\underline{\a\1(a\a\1(a')_{R})_{2}}, a''_{1})\\
 &&~\times \upsilon(\underline{\b\1(\b\1(b_{R})b')_{1}}, b''_{2})\psi(\underline{\b\1(\b\1(b_{R})b')_{2}}, a''_{2})\\
 &\stackrel{(HA1)(HC1)}{=}&\varphi(\a\1(a_{1})\a\1(\a\1(a')_{R1}), b''_{1})
 \tau(\a\1(a_{2})\a\1(\a\1(a')_{R2}), a''_{1})\\
 &&~\times \upsilon(\b^{-2}(b_{R1})\b\1(b'_{1}), b''_{2})
 \psi(\b^{-2}(b_{R2})\b\1(b'_{2}), a''_{2})\\
 &\stackrel{(CHA2)(SP2)}{=}&\varphi(a_{1}, \b\1(b''_{11}))\varphi(\underline{\a\1(a')_{R1}}, \b\1(b''_{12}))
 \tau(a_{2}, \a\1(a''_{11}))\\
 &&\times \tau(\underline{\a\1(a')_{R2}}, \a\1(a''_{12}))\upsilon(\b\1(\underline{b_{R1}}), \b\1(b''_{21}))\upsilon(b'_{1}, \b\1(b''_{22}))\\
   \end{eqnarray*}
  \begin{eqnarray*}
 &&\times \psi(\b\1(\underline{b_{R2}}), \a\1(a''_{21}))\psi(b'_{2}, \a\1(a''_{22}))\\
 &\stackrel{}{=}&\varphi(a_{1}, \underline{\b\1(b''_{11})})\varphi(\a\1(a')_{1R}, \underline{\b\1(b''_{12})})
 \tau(a_{2}, \underline{\a\1(a''_{11})})\\
 &&\times \tau(\a\1(a')_{2r}, \underline{\a\1(a''_{12})})\times \upsilon(\b\1(b_{1R}), \underline{\b\1(b''_{21})})\upsilon(b'_{1}, \underline{\b\1(b''_{22})})\\
 &&\times \psi(\b\1(b_{2r}), \underline{\a\1(a''_{21})})\psi(b'_{2}, \underline{\a\1(a''_{22})})\\
 &\stackrel{(3)}{=}&\varphi(a_{1}, b''_{1})\varphi(\underline{\a\1(a')_{1R}}, \underline{\b^{-2}(b''_{211})})
 \tau(a_{2}, a''_{1})\\
 &&\times \tau(\underline{\a\1(a')_{2r}}, \underline{\a^{-2}(a''_{211})})\upsilon(\underline{\b\1(b_{1R})}, \underline{\b^{-2}(b''_{212})})\upsilon(b'_{1}, \b\1(b''_{22}))\\
 &&\times \psi(\underline{\b\1(b_{2r})}, \underline{\a^{-2}(a''_{212})})\psi(b'_{2}, \a\1(a''_{22}))\\
 &\stackrel{(4)(HC1)}{=}&\varphi(a_{1}, b''_{1})\underline{\varphi(\a(\a\1(\a\1(a')_{1})_{R}), \b^{-2}(b''_{21})_{1})}
 \tau(a_{2}, a''_{1})\\
 &&\times \underline{\tau(\a(\a\1(\a\1(a')_{2})_{r}), \a^{-2}(a''_{21})_{1})}
 ~\underline{\upsilon(\b\1(b_{1})_{R}, \b^{-2}(b''_{21})_{2})}\\
 &&\times \upsilon(b'_{1}, \b\1(b''_{22}))\underline{\psi(\b\1(b_{2})_{r}, \a^{-2}(a''_{21})_{2})}\psi(b'_{2}, \a\1(a''_{22}))\\
 &\stackrel{(D1)(D2)}{=}&\varphi(a_{1}, \underline{b''_{1}})\varphi(a'_{1}, \underline{\b^{-2}(b''_{21})_{2}})
 \tau(a_{2}, \underline{a''_{1}})\tau(a'_{2}, \underline{\a^{-2}(a''_{21})_{2}})\\
 &&\times \upsilon(b_{1}, \underline{\b^{-2}(b''_{21})_{1}})\upsilon(b'_{1}, \underline{\b\1(b''_{22})})
 \psi(b_{2}, \underline{\a^{-2}(a''_{21})_{1}})\psi(b'_{2}, \underline{\a\1(a''_{22})})\\
 &\stackrel{(3)}{=}&\varphi(a_{1}, \underline{\b\1(b''_{11})})\varphi(a'_{1}, \underline{\b\1(b''_{21})})
 \tau(a_{2}, \underline{\a\1(a''_{11})})\tau(a'_{2}, \underline{\a\1(a''_{21})})\\
 &&\times \upsilon(b_{1}, \underline{\b\1(b''_{12})})\upsilon(b'_{1}, \underline{\b\1(b''_{22})})
 \psi(b_{2}, \underline{\a\1(a''_{12})})\psi(b'_{2}, \underline{\a\1(a''_{22})})\\
 &\stackrel{(HC1)}{=}&\varphi(a_{1}, \b\1(b''_{1})_{1})\tau(a_{2}, \a\1(a''_{1})_{1})
 \upsilon(b'_{1}, \b\1(b''_{2})_{2})\psi(b'_{2}, \a\1(a''_{2})_{2})\\
 &&\times \varphi(a'_{1}, \b\1(b''_{2})_{1})\tau(a'_{2}, \a\1(a''_{2})_{1})
 \upsilon(b_{1}, \b\1(b''_{1})_{2})\psi(b_{2}, \a\1(a''_{1})_{2})\\
 &\stackrel{}{=}&\s(\a(a)\o \b(b), a''_1\o b''_1)\s(\a(a')\o \b(b'), a''_2\o
 b''_2).
 \end{eqnarray*}

 $(CHA3)$ for $\s$ can be proved by similar method.

 And we check $(CHA4)$ as follows. For all $a,a'\in A$ and $b,b'\in B$,
 \begin{eqnarray*}
 &&\s(a_{1}\o b_{1}, a'_{1}\o b'_{1})(a_{2}\o b_{2})(a'_{2}\o b'_{2})\\
 &\stackrel{}{=}&u
 (\a\1(a_{1})_{1}, \b\1(b'_{1})_{1})\tau(\a\1(a_{1})_{2}, \a\1(a'_{1})_{1})
 \upsilon(\b\1(b_{1})_{1}, \b\1(b'_{1})_{2})\\
 &&\times \psi(\b\1(b_{1})_{2}, \a\1(a'_{1})_{2})(a_{2}\a\1(a'_{2})_{R}\o \b\1(b_{2R})b'_{2})\\
 &\stackrel{(HC1)}{=}&\varphi(\a\1(a_{11}), \b\1(b'_{11}))\tau(\a\1(a_{12}), \a\1(a'_{11}))
 \upsilon(\b\1(b_{11}), \b\1(b'_{12}))\\
 &&\times \psi(\b\1(b_{12}), \a\1(a'_{12}))(a_{2}\a\1(a'_{2})_{R}\o \b\1(b_{2R})b'_{2})\\
 &&\times \psi(\b\1(b_{1})_{2}, \a\1(a'_{1})_{2})(a_{2}\a\1(a'_{2})_{R}\o \b\1(b_{2R})b'_{2})\\
 &\stackrel{(2)}{=}&\varphi(a_{1}, b'_{1})\tau(\a\1(a_{21}), a'_{1})
 \upsilon(b_{1}, \b\1(b'_{21}))\psi(\b\1(b_{21}), \a\1(a'_{21}))\\
 &&\times  (\a\1(a_{22})\a^{-2}(a'_{22})_{R}\o \b\1(\b\1(b_{22})_{R})\b\1(b'_{22}))\\
 &\stackrel{(HC1)}{=}&\varphi(a_{1}, b'_{1})\tau(\a\1(a_{2})_{1}, a'_{1})
 \upsilon(b_{1}, \b\1(b'_{2})_{1})\underline{\psi(\b\1(b_{2})_{1}, \a\1(a'_{2})_{1})}\\
 &&\times  (\a\1(a_{2})_{2}\underline{\a\1(\a\1(a'_{2})_{2})_{R}}\o
 \underline{\b\1(\b\1(b_{2})_{2R})}\b\1(b'_{2})_{2})\\
  \end{eqnarray*}
  \begin{eqnarray*}
 &\stackrel{(D5)}{=}&\varphi(a_{1}, b'_{1})\tau(\a\1(a_{2})_{1}, a'_{1})
 \upsilon(b_{1}, \b\1(b'_{2})_{1})\psi(\b\1(b_{2})_{2}, \a\1(a'_{2})_{2})\\
 &&\times  (\a\1(a_{2})_{2}\a\1(a'_{2})_{1}\o
 \b\1(b_{2})_{1}\b\1(b'_{2})_{2})\\
 &\stackrel{(1)(HC1)}{=}&\varphi(a_{1}, b'_{1})\underline{\tau(\a\1(a_{2})_{1}, \a\1(a'_{1})_{1})}
 ~\underline{\upsilon(\b\1(b_{1})_{1}, \b\1(b'_{2})_{1})}\psi(b_{2}, a'_{2})\\
 &&\times  (\underline{\a\1(a_{2})_{2}\a\1(a'_{1})_{2}}\o
 \underline{\b\1(b_{1})_{2}\b\1(b'_{2})_{2}})\\
 &\stackrel{(CHA4)}{=}&u
 (a_{1}, b'_{1})\tau(\a\1(a_{2})_{2}, \a\1(a'_{1})_{2})
 \upsilon(\b\1(b_{1})_{2}, \b\1(b'_{2})_{2})\psi(b_{2}, a'_{2})\\
 &&\times  (\a\1(a'_{1})_{1}\a\1(a_{2})_{1}\o
 \b\1(b'_{2})_{1}\b\1(b_{1})_{1})\\
 &\stackrel{(1)(HC1)}{=}&\underline{\varphi(\a\1(a_{1})_{1}, \b\1(b'_{1})_{1})}
 \tau(a_{2}, \a\1(a'_{1})_{2})
 \upsilon(\b\1(b_{1})_{2}, b'_{2})\psi(b_{2}, a'_{2})\\
 &&\times  (\a\1(a'_{1})_{1}\underline{\a\1(a_{1})_{2}}\o
 \underline{\b\1(b'_{1})_{2}}\b\1(b_{1})_{1})\\
 &\stackrel{(D6)}{=}&\varphi(\a\1(a_{1})_{2}, \b\1(b'_{1})_{2})
 \tau(a_{2}, \a\1(a'_{1})_{2})
 \upsilon(\b\1(b_{1})_{2}, b'_{2})\psi(b_{2}, a'_{2})\\
 &&\times  (\a\1(a'_{1})_{1}\a\1(\a\1(a_{1})_{1})_R\o
 \b\1(\b\1(b'_{1})_{1R})\b\1(b_{1})_{1})\\
 &\stackrel{(2)(3)}{=}&(a'_{1}\a\1(a_{1})_R\o
 \b\1(b'_{1R})b_{1})\varphi(\a\1(a_{2})_{1}, \b\1(b'_{1})_{2})\\
 &&\times  \tau(\a\1(a_{2})_{2}, \a\1(a'_{2})_{1})\upsilon(\b\1(b_{2})_{1}, \b\1(b'_{2})_{2})\psi(\b\1(b_{2})_{2},
 \a\1(a'_{2})_{2})\\
 &\stackrel{}{=}&(a'_{1}\o b'_{1})(a_{1}\o b_{1})\s(a_{2}\o b_{2}, a'_{2}\o
 b'_{2}).
 \end{eqnarray*}

 Therefore, $(A\natural_R B,\a\o \b,\s)$ is a cobraided Hom-Hopf algebra.
 \hfill $\square$

\smallskip

 Thus it follows from Lemmas 4.2-4.7 that we have

 {\bf Theorem 4.8} $R$-smash product Hom-Hopf algebra $(A\natural_R B,\a\o \b)$ is cobraided
 if and only if there exist forms $\tau: A\o A\lr K$, $\varphi: A\o B\lr K$, $\psi: B\o A\lr K$,
 and $\upsilon: B\o B\lr K$ such that $(A,\a,\tau)$ and $(B,\b,\upsilon)$ are cobraided Hom-Hopf algebras,
 $(A,B,\varphi)$ and $(B,A,\psi)$ are Hom-skew pairings, the conditions $(D1)-(D6)$ in Lemma 4.5 hold.
 Moreover, the cobraided structure $\s$ on $(A\natural_R B,\a\o \b)$ has a decomposition
 $$
 \s(\a(a)\o \b(b), \a(a')\o \b(b'))=\varphi(a_1, b'_1)\tau(a_2, a'_1)\upsilon(b_1, b'_2)\psi(b_2,
 a'_2).
 $$

\smallskip

 {\bf Theorem 4.9} Smash product Hom-Hopf algebra $(A\natural H,\a\o \b)$ is cobraided
 if and only if there exist forms $\tau: A\o A\lr K$, $\varphi: A\o H\lr K$, $\psi: H\o A\lr K$,
 and $\upsilon: H\o H\lr K$ such that $(A,\a,\tau)$ and $(H,\b,\upsilon)$ are cobraided Hom-Hopf algebras,
 $(A,H,\varphi)$ and $(H,A,\psi)$ are Hom-skew pairings, the conditions $(D1)'-(D6)'$ below hold.
 For all $a,a'\in A$ and $h,h'\in B$,
 \begin{eqnarray*}
 &(D1)'& \varphi(\b(h'_1)\rhd a, h_1)\upsilon(h'_2, h_2)=\upsilon(\b(h'), h_1)\varphi(\a(a), h_2),\\
 &(D2)'& \tau(\b(h_1)\rhd a, a'_1)\psi(h_2, a'_2)=\psi(\b(h), a'_1)\tau(\a(a), a'_2),\\
 &(D3)'& \upsilon(h_1, h'_2)\psi(h_2, \b(h'_1)\rhd a)=\psi(h_1, \a(a)\upsilon(h_2, \b(h')),\\
 &(D4)'& \varphi(a_1, h_2)\tau(a_2, \b(h_1)\rhd a')=\tau(a_1, \a(a')\varphi(a_2, \b(h)),\\
   \end{eqnarray*}
  \begin{eqnarray*}
 &(D5)'& \psi(h_1, a_1)(\b(h_{21})\rhd a_2\o h_{22})=(\a(a_1)\o \b(h_{1}))\psi(h_2, a_2),\\
 &(D6)'& \varphi(a_1, h_1)(\a(a_2)\o \b(h_{2}))=(\b(h_{11})\rhd a_1\o h_{12})\varphi(a_2, h_2).
 \end{eqnarray*}
 Moreover, the cobraided structure $\s'$ on $(A\natural H,\a\o \b)$ has a decomposition
 $$
 \s'(\a(a)\o \b(h), \a(a')\o \b(h'))=\varphi(a_1, h'_1)\tau(a_2, a'_1)\upsilon(h_1, h'_2)\psi(h_2,
 a'_2).
 $$

 {\bf Proof} Let $R(h\o a)=h_{1}\rhd a\o h_{2}, \forall a\in A, h\in H$ in Theorem 4.8. \hfill $\square$

\section{Applications}

 In this section, we give the applications of the main results in Sec.3 and 4 to a concrete example.

 The following result is clear.

 {\bf Lemma 5.1} Let $K\mathbb{Z}_2=K\{1,a\}$ be Hopf group algebra (see \cite{Sw}). Then $(K\mathbb{Z}_2, id_{K\mathbb{Z}_2}, \upsilon)$ is a cobraided Hom-Hopf algebra, where $\upsilon$ is given by
 \begin{center}
 \begin{tabular}{c|cc}
  $\upsilon$ & 1 & a  \\
  \hline
  1& 1 & 1 \\
  a & 1& $-1$\\
 \end{tabular}~.
 \end{center}

\smallskip

 Let $T_{2,-1}=K\{1, g, x, gx|g^2=1, x^2=0, xg=-gx\}$ be Taft's Hopf algebra (see \cite{Ta}), its coalgebra structure and antipode are given by
 $$
 \D(g)=g\o g,~\D(x) = x\o g+1\o x, ~\D(gx)=gx\o 1+g\o gx;
 $$
 $$
 \v(g) = 1, \v(x) = 0, \v(gx)=0;
 $$
 and
 $$
 S(g)=g,~S(x)=gx,~S(gx)=-x.
 $$

 Define a linear map $\a$: $T_{2,-1}\lr T_{2,-1}$ by
 $$
 \a(1)=1,~\a(g)=g,~\a(x)=kx,~\a(gx)=kgx
 $$
 where $0\neq k\in K$. Then $\a$ is an automorphism of Hopf algebras.

 So we can get a Hom-Hopf algebra $H_{\a}=(T_{2,-1}, \a\ci \mu_{T_{2,-1}}, 1_{T_{2,-1}}, \D_{T_{2,-1}}\ci \a, \v_{T_{2,-1}}, \a)$ (see \cite{MS2}).

  {\bf Lemma 5.2} Let $H_{\a}$ be the Hom-Hopf algebra defined as above. Then $(H_{\a}, \a, \tau)$ is a cobraided Hom-Hopf algebra, where $\tau$ is given by
 \begin{center}
 \begin{tabular}{c|cccc}
  $\tau$ & 1 & g  & x & gx\\
  \hline
  1& 1 & 1 & 0 & 0\\
  g & 1& $-1$& 0 & 0\\
  x & 0 & 0 & 0 & 0\\
  gx & 0&0& 0 & 0\\
 \end{tabular}~.
 \end{center}

 {\bf Proof} It is straightforward by a tedious computation.       \hfill $\square$

\smallskip

 {\bf Theorem 5.3} Let $K\mathbb{Z}_2$ be Hopf group algebra and $H_{\a}$ be the Hom-Hopf algebra defined as above.
 Define module action $\rhd: K\mathbb{Z}_2\o H_{\a} \lr H_{\a}$ by
 $$
 1_{K\mathbb{Z}_2}\rhd 1_{H_{\a}}=1_{H_{\a}},~1_{K\mathbb{Z}_2}\rhd g=g,
 $$
  $$
 1_{K\mathbb{Z}_2}\rhd x=kx,~1_{K\mathbb{Z}_2}\rhd gx=kgx,
 $$
 $$
 a\rhd 1_{H_{\a}}=1_{H_{\a}},~a\rhd g=g,
 $$
 $$
 a\rhd x=-kx,~a\rhd gx=-kgx,
 $$
 Then by a routine computation we can get $H_{\a}$ is a $K\mathbb{Z}_2$-module Hom-algebra. Therefore, by Theorem 3.3,
 $(H_{\a}\natural K\mathbb{Z}_2, \a\o id_{K\mathbb{Z}_2})$ is a smash product Hom-algebra.

 Furthermore, $(H_{\a}\natural K\mathbb{Z}_2, \a\o id_{K\mathbb{Z}_2})$ with the tensor product Hom-coalgebra becomes a Hom-Hopf algebra, where the antipode $\bar{S}$ is given by
 $$
 \bar{S}(1_{H_{\a}}\o 1_{K\mathbb{Z}_2})=1_{H_{\a}}\o 1_{K\mathbb{Z}_2},~\bar{S}(1_{H_{\a}}\o a)=1_{H_{\a}}\o a,
 $$
 $$
 \bar{S}(g\o 1_{K\mathbb{Z}_2})=g\o 1_{K\mathbb{Z}_2},~\bar{S}(g\o a)=g\o a,
 $$
 $$
 \bar{S}(x\o 1_{K\mathbb{Z}_2})=-gx\o 1_{K\mathbb{Z}_2},~\bar{S}(x\o a)=-gx\o a,
 $$
 $$
 \bar{S}(gx\o 1_{K\mathbb{Z}_2})=x\o 1_{K\mathbb{Z}_2},~\bar{S}(gx\o a)=x\o a.
 $$

\smallskip

 {\bf Lemma 5.4} Let $K\mathbb{Z}_2$ be the Hopf group algebra and $H_{\a}$  be the Hom-Hopf algebra defined as above.
 Define two linear maps $\varphi: H_{\a}\o K\mathbb{Z}_2\lr K$ and $\psi: K\mathbb{Z}_2\o H_{\a}\lr K$ as follows:
 \begin{center}
 \begin{tabular}{c|cc}
  $\varphi$ & 1 & a  \\
  \hline
  1& 1 & 1 \\
  g & 1& $-1$ \\
  x& 0 & 0 \\
  gx & 0& 0 \\
 \end{tabular}
~~~~~~~~~~
  \begin{tabular}{c|cccc}
  $\psi$ & 1 & g  & x & gx\\
  \hline
  1& 1 & 1& 0 & 0  \\
  a & 1& $-1$& 0 & 0  \\
 \end{tabular}.
 \end{center}
 Then $(H_{\a},K\mathbb{Z}_2,\varphi)$ and $(K\mathbb{Z}_2,H_{\a},\psi)$ are two Hom-skew pairings.

 {\bf Proof} Straightforward.       \hfill $\square$

\smallskip

 {\bf Theorem 5.5} With notations as above. Smash product Hom-Hopf algebra $(H_{\a}\natural K\mathbb{Z}_2, \a\o id_{K\mathbb{Z}_2},\s)$ is a cobraided Hom-Hopf algebra with the cobraiding $\s$ is given as follows:
 \begin{center}
 \begin{tabular}{c|cccccccc}
  $\s$ & $1\o 1$ & $1\o a$  & $g\o 1$ &$g\o a$&$ x\o 1$ & $x\o a$  & $gx\o 1$ & $gx\o a$\\
  \hline
  $1\o 1$& 1 & 1& 1 & 1 & 0 & 0& 0 & 0 \\
  $1\o a$ & 1& $-1$& $-1$ & 1 & 0 & 0& 0 & 0 \\
  $g\o 1$ & 1& $-1$& $-1$ & 1  & 0 & 0& 0 & 0\\
  $g\o a$ & 1 & 1& 1 & 1& 0 & 0& 0 & 0 \\
  $x\o 1$& 0 & 0& 0 & 0 & 0 & 0& 0 & 0 \\
  $x\o a$ & 0 & 0& 0 & 0 & 0 & 0& 0 & 0 \\
  $gx\o 1$& 0 & 0& 0 & 0  & 0 & 0& 0 & 0\\
  $gx\o a$ & 0 & 0& 0 & 0 & 0 & 0& 0 & 0 \\
 \end{tabular}.
 \end{center}

 {\bf Proof} It is easy to prove that the conditions $(D1)'-(D6)'$ hold. And by Lemma 5.1,5.2,5.4 and Theorem 4.9, we can finish the proof.       \hfill $\square$

\subsection*{Acknowledgements}
 The authors are deeply indebted to the referee for his/her very useful suggestions and some improvements to the original manuscript, and for bringing our attention to ref\cite{CWZ}. This work was partially supported by the NNSF of China (No. 11101128).

\end{document}